\documentclass[a4paper]{amsart}
\pdfoutput=1
\usepackage{amsthm}
\usepackage{amssymb}
\usepackage{bm}
\usepackage{bbding}
\usepackage{biblatex}
	\bibliography{bibliography.bib}

\theoremstyle{definition}

\theoremstyle{plain}

\newtheorem{cor}{Corrolary}

\title[Identities for $c_{mk}$- and $a_{mk}$-Weighted Sums]{Further Identities for $c_{mk}$- and $a_{mk}$-Weighted Sums and a Remark on a Representation of Pythagoras' Equation}

\begin{document}
\author[C. Muschielok]{Christoph Muschielok}
\address{Department of Chemistry, 
         Technical University of Munich, 
         Lichtenbergstra{\ss}e 4,
         85748 Garching, 
         Germany}
\email{c.muschielok(at)tum.de}

\begin{abstract}
	We present some properties of the expansion coefficients~$a_{mk}$ and
$c_{mk}$ of a pair of dual bases, 
	\[ n^m = \sum_{k=2}^m c_{mk} \psi_{k}(n)\,, \] 
and 
	\[ \psi_m(n) = n + (m-1)(n-1) B_{n-1,m-1}\,, \]
we introduced earlier in arXiv:2207.01935v1. Here, $B_{a,b} = (a+b)!/(a!\,b!)$
is a binomial coefficient. We extend the knowledge on the $c_{mk}$ by giving an
explicit expression for them in terms of the Stirling numbers of
the second kind.
From the interchangeability of the indices of the binomial coefficients, we
formulate the central identity we use in this work:
	\[ \psi_m(n) - n = \psi_n(m) - m\,. \]
With this equation, we evaluate sums of the form 
	\[ T^\alpha_m = \sum_{k=2}^m c_{mk} k^\alpha\,.\]
Explicitly, the case~$T_m^1$ is handled. Furthermore, we indicate connections 
of~$T_m^2$ and~$T_m^3$ to the Mersenne numbers (general integer exponent) and 
the OEIS entry A024023. We conclude with a small remark on how we can represent
Pythagoras' equation in terms of the $a_{mk}$ coefficients.
\end{abstract}
\maketitle

\section{Introduction}
Prior\cite{muschielok2022powersums}, we have shown that a power sum of $a$-th
order (or hyper-sum),
\begin{align}
	S^{(a)}_m(n) = \sum_{\nu_a = 1}^n \cdots 
                       \sum_{\nu_2 = 1}^{\nu_3} 
                       \sum_{\nu_1=1}^{\nu_2} \nu_1^m\,,
\label{eq:PowerSum}
\end{align}
can be replaced by a polynomial expansion
\begin{align}
	S^{(a)}_m(n) = \sum_{k=2}^m c_{mk}\psi^{(a)}_k(n)\,,
\end{align}
where the polynomials~$\psi^{(a)}_m(n)$ are given by
\begin{align}
	\psi_m^{(a)}(n) = B_{a+1,n-1} 
                          + \frac{m(m-1)}{m+a} (n-1) B_{m+a-1,n-1}\,.
\end{align}
The expansion coefficients~$c_{mk}$\cite{oeis2022cmk} and polynomials are
introduced and discussed in detail in Ref.~\cite{muschielok2022powersums}. 
$B_{ab} = (a+b)!/(a!\,b!)$ denotes a binomial coefficient.

In the following, if
not stated otherwise, we require~$n\in\mathbb{N}$ and~$m \ge 2$ (although we
might readily include~$m=1$, this case is considered separately as stated in
Ref.~\cite{muschielok2022powersums}).

In the following we furthermore consider only the $\psi$-polynomials for~$a=0$,
\begin{align}
	\psi_m^{(0)}(n) = n + (m-1)(n-1)B_{n-1,m-1} =: \psi_m(n)\,.
	\label{eq:ExplicitPsi}
\end{align}
Trivially, we can expand the polynomials~$\psi_m(n)$ in powers of~$n$. Through
this expansion, we obtain the coefficients~$a_{mk}$:
\begin{align}
	\psi_m(n) = \sum_{k=2}^m a_{mk} n^k\,,
	\label{eq:nExpansion}
\end{align}
where we have already argued using Vieta's formulas, that the
coefficients~$a_{m0}$ and~$a_{m1}$ are both equal to
zero.\cite{muschielok2022powersums} As J.~L.~Cereceda points
out in a nice addendum to the author's work,\cite{cereceda2022powersums} they 
have the closed form
\begin{align}
	a_{mk} = \frac{{m-1 \brack k-1} - {m-1 \brack k}}{(m-2)!}\,,
\end{align}
where $n \brack m$ is an unsigned Stirling number of the first kind. Using the
recursion relation for the Stirling numbers of the first kind (cf. Table~250 in
Ref.~\cite{graham1994concrete}),
\begin{align}
	{n \brack k} = {n - 1 \brack k} (n-1) + {n-1 \brack k-1}\,,
\end{align}
we may write down a similar relation for the
coefficients~$a_{mk}$,
\begin{align}
	a_{mk} = a_{m-1, k} + \frac{a_{m-1,k-1}}{m-2}\,.
\end{align}

The dual expansion, i.e. the expansion of the monomial~$n^m$ in terms of
the~$\psi_m(n)$,
\begin{align}
	n^m = \sum_{k=2}^m c_{mk} \psi_k(n)\,,
	\label{eq:PsiExpansion}
\end{align}
defines the coefficients~$c_{mk}$ which can be used in the power sum in
Eq.~\eqref{eq:PowerSum}. 

Cereceda gives furthermore several identities and properties which the author 
had not touched yet in his published work. In particular, here, we will use the 
following two identities given by Cereceda (Eqs. (4) and (8) \emph{ibid.}):
\begin{align}
	\sum_{k=2}^m c_{mk} = 1 \label{eq:TrivialSumCmk}
\end{align}
and
\begin{align}
	\sum_{k=2}^m c_{mk} B_{k-1,n-1} = \frac{n^{m-1} - n}{n-1}\,.
	\label{eq:WeightedBinomialCoeff}
\end{align}
Those are easily obtained by inserting $n=1$ into Eq.~\eqref{eq:PsiExpansion}
and by inserting the explicit form of $\psi_m(n)$, Eq.~\eqref{eq:ExplicitPsi},
respectively.

The equivalent identity to Eq.~\eqref{eq:TrivialSumCmk} is given by inserting
$n=1$ into Eq.~\eqref{eq:nExpansion}:
\begin{align}
	\sum_{k=2}^m a_{mk} = 1\,. \label{eq:TrivialSumAmk}
\end{align}
From Eqs.~\eqref{eq:TrivialSumCmk} and~\eqref{eq:TrivialSumAmk} it is clear,
that any constant term~$C$ may be pulled inside a sum of sequence
elements~$A_k$ which is weighted by~$c_{mk}$ or~$a_{mk}$, that is
\begin{align}
	\sum_{k=2}^m c_{mk} A_k - C = \sum_{k=2}^m c_{mk} (A_k - C)
\intertext{and}
	\sum_{k=2}^m a_{mk} A_k - C = \sum_{k=2}^m a_{mk} (A_k - C)\,.
\end{align}

We can read the two expansions, Eqs.~\eqref{eq:nExpansion}
and~\eqref{eq:PsiExpansion}, as vector--matrix products of a coefficient vector
and a vector of polynomial values. The latter vector can be interpreted as the
left-product of what one might call here a ``vector of choice'' $(0, 0, \dotsc,
1, 0, \dotsc, 0)^T$ and a polynomial matrix. We denote the coefficient vectors
by $c_m = (c_{m\mu})_{2\le \mu \le m}$ and $a_m = (a_{m\mu})_{2\le \mu \le m}$.
For the matrices we write $\Psi_m(n) = (\psi_\mu(\nu))_{2\le\mu\le m,
1\le\nu\le n}$ and $V_m(n) = (n^m)_{2\le\mu\le m, 1\le\nu\le n}$ (as this is a
kind of Vandermonde matrix). With this, we write for the right-product of the
polynomial matrix with the coefficient vector
\begin{align}
	&(1, 2^m, \dotsc, n^m) = c_m^T \Psi_m(n)
\intertext{and}
	&(1, \psi_m(2), \dotsc, \psi_m(n)) = a_m^T V_m(n)\,.
\end{align}
In the following, we consider the triple products or matrix elements ($\alpha, m \ge 2$)
\begin{align}
	T_m^\alpha := c_m^T \Psi_m(\alpha) c_\alpha 
                     = \sum_{l=2}^\alpha c_{\alpha l} l^\alpha\,.
\end{align}
The corresponding matrix elements
\begin{align}
	U_m^\alpha := a_m^T V_m(\alpha) a_\alpha 
                    = \sum_{l=2}^\alpha a_{\alpha l} \psi_m(l),\quad(\alpha, m \ge 2)\,,
\end{align}
will not be discussed further in this manuscript. On the other hand, we will
explicitly handle the special cases~$T_m^1$ and~$U_m^1$.

\section{Explicit form for the $c_{mk}$ coefficients}
We start from Eq.~\eqref{eq:WeightedBinomialCoeff}, which we want to write in
the following way:
\begin{align}
	\sum_{k=2}^m c_{mk} \frac{n^{\overline{k-1}}}{(k-2)!} = \frac{n^m - n}{n-1}\,,
	\label{eq:RisingFactorialExpansion}
\end{align}
where~$n^{\overline{k-1}} = \prod_{j=0}^{k-2} (n+j)$ means Knuth's notation for
the rising factorial.

We can expand the right-hand side using polynomial division:
\begin{align}
	\frac{n^m - n}{n-1} = n^{m-1} + n^{m-2} + \dotsc + n = \sum_{l=2}^{m} n^{l-1}\,.
	\label{eq:PolynomialExpansion}
\end{align}
We know that we can express powers in terms of the falling or rising factorials
by means of the Stirling numbers of the second kind. Using the equation\cite[p.
250]{graham1994concrete}
\begin{align}
	x^n = \sum_k {n \brace k} (-1)^{n-k} x^{\overline{k}}\,,
\end{align}
to substitute the powers in Eq.~\eqref{eq:PolynomialExpansion} and replacing the
right-hand side in Eq.~\eqref{eq:RisingFactorialExpansion}, we have
\begin{align}
	\sum_{k=2}^m c_{mk} \frac{n^{\overline{k-1}}}{(k-2)!} = \sum_{l=2}^m \left(\sum_j {l-1 \brace j} (-1)^{l-1-j} \right) n^{\overline{j}}\,.
\end{align}
We switch the summation order and resolve a factor of 1 into $(j-1)!/(j-1)!$:
\begin{align}
	\sum_{k=2}^m c_{mk} \frac{n^{\overline{k-1}}}{(k-2)!} &= \sum_j \left(\sum_{l=2}^m {l-1 \brace j} (-1)^{l-1-j} (j-1)! \right) \frac{n^{\overline{j}}}{(j-1)!}\,.
\end{align}
By comparison of coefficients ($k-1=j$), we find
\begin{align}
	c_{mk} =  (k-2)!\sum_{l=2}^m {l-1 \brace k-1} (-1)^{l-k}\,.
\end{align}

\section{Exchanging $m$ and $n$ in $\psi_m(n)$}
Before we present the identities for the sums shown in the introductory section, 
we want to consider a property of the polynomials~$\psi_m(n)$, which deems
important for deriving these identities, first.
\begin{cor}
	\label{cor:Commutator}
	The polynomial~$\psi_m(n) - n$ is symmetric under exchange of $n$ and
	$m$, that is
	\begin{align}
		\psi_m(n) - n = \psi_n(m) - m
	\intertext{or}
		\psi_m(n) = \psi_n(m) + n - m\,.
	\end{align}
\end{cor}
\emph{Nota bene}: we can see that through this, we are able to extend the
definition of the function~$\psi_m(n)$ to $m=0$.
\begin{align}
	\psi_0(n) := \lim_{m \rightarrow 0} (\psi_n(m) - m + n) 
                   = n + \lim_{m\rightarrow 0} \psi_n(m)\,,
\end{align}
where we can find the remaining limit by considering the value of the
polynomial~$\psi_m(n)$, in contrast to its binomial coefficient form (which is
undefined for $m=0$). If the limit of the polynomial ($n\in\mathbb{R}$) exists,
so does the limit of the sequence. Clearly, the limit is zero. Thus, we find 
that~$\psi_0(n) = \psi_1(n)$, that is
\begin{align}
	\psi_0(n) = n\,.
\end{align}

\section{The Weighted Sums}
Next we want to present the values of the sums
\begin{align*}
	T_m^1 = \sum_{k=2}^m c_{mk} k
\quad\text{and}\quad
	U_m^1 = \sum_{l=2}^m a_{ml} l\,.
\end{align*}
First, we want to write Eq.~\eqref{eq:WeightedBinomialCoeff} in a form which
stresses that we are dealing with a polynomial. We begin with
Eq.~\eqref{eq:PsiExpansion} and insert the explicit form for~$\psi_m(n)$:
\begin{align}
	\sum_{k=2}^m c_{mk} \left[ n 
                                   + (n-1) 
                                     \frac{1}{(k-2)!} 
                                     \prod_{l=0}^{k-2} (n + l) 
                           \right] 
                           = n^m\,.
\end{align}
After further manipulation, this equation becomes
\begin{align}
	\sum_{k=2}^m c_{mk} \frac{1}{(k-2)!} \prod_{l=0}^{k-2} (n+l) 
		= \frac{n^m - n}{n - 1}\,.
\end{align}
Next, consider the limit for~$n\rightarrow 1$:
\begin{align}
	\sum_{k=2}^m c_{mk} \frac{1}{(k-2)!} \prod_{l=0}^{k-2} (l+1) 
		= \lim_{n \rightarrow 1} \frac{n^m - n}{n-1}\,.
\end{align}
On the left-hand side, we can further cancel the product term with the
factorial and arrive at
\begin{align}
	\sum_{k=2}^m c_{mk} (k - 1) = \lim_{n \rightarrow 1} \frac{n^m - n}{n-1}\,.
\end{align}
By means of Eq.~\eqref{eq:TrivialSumCmk}, this can be further simplified to
\begin{align}
	\sum_{k=2}^m c_{mk} k &= 1 + \lim_{n \rightarrow 1} \frac{n^m - n}{n-1}\,.
\end{align}
If the limit exists for the rational function $(x^m - x) / (x-1)$,
$x\in\mathbb{R}$, its value will be the same for the sequence we are dealing
with, originally. We have
\begin{align}
	\lim_{n\rightarrow 1} \frac{n^m - n}{n-1} = 
	    \lim_{x \rightarrow 1} \frac{x^m - x}{x - 1} = 
            \lim_{x \rightarrow 1} (m x^{m-1} - 1)\,,
\end{align}
where in the last step we used l'H\^{o}pital's rule. Finally, we have the
equation
\begin{align}
	\sum_{k=2}^m c_{mk} k = m\,. \label{eq:WeightedSumCmk}
\end{align}

For finding the corresponding expression in which the~$a_{mk}$ coefficients are
used, we can exploit the orthogonality of the coefficient vectors,
\begin{align}
	\sum_{k=2}^m a_{ml} c_{lk} = \delta_{mk}\,. \label{eq:Ortho}
\end{align}
Multiplying with $a_{mk}$ from the left on Eq.~\eqref{eq:WeightedSumCmk} and 
summing over~$k$, we have
\begin{align}
	\sum_{k=2}^m a_{mk} k &= \sum_{k=2}^m a_{mk} \sum_{l=2}^k c_{kl} l\,.
\end{align}
We can decouple the summation boundaries by using~$c_{kl} = 0$ for 
$l < 2 \lor l > m$:\cite{muschielok2022powersums}
\begin{align}
	\sum_{k=2}^m a_{mk} k = \sum_{k=2}^m a_{mk} \sum_{l=2}^m c_{kl} l\,.
\end{align}
Then, after interchanging the summation order and using Eq.~\eqref{eq:Ortho},
we finally have
\begin{align}
	\sum_{k=2}^m a_{mk} k = m\,.
\end{align}

\section{The matrix elements $T_m^\alpha$ and $U_m^\alpha$ for $\alpha > 1$}
At last we consider~$T_m^\alpha$, for $\alpha > 1$. Again, we use
Eq.~\eqref{eq:PsiExpansion} to expand~$k^\alpha$:
\begin{align}
	\sum_{k=2}^m c_{mk} k^\alpha = \sum_{k=2}^m c_{mk} 
                                       \sum_{l=2}^\alpha c_{\alpha l} \psi_l(k)\,.
\end{align}
Substract and add~$k$ within the sum over~$l$. By means of
Corrolary~\ref{cor:Commutator}, we interchange~$k$ and~$l$ for~$\psi_{l}(k)-k$:
\begin{align}
	\sum_{k=2}^m c_{mk} 
        \sum_{l=2}^\alpha c_{\alpha l} [\psi_l(k) -k +k]
	    = \sum_{l=2}^\alpha 
              \sum_{k=2}^m c_{\alpha l} c_{mk} [\psi_k(l) - l + k]\,.
\end{align}
Using again the expansion in terms of the $c_{mk}$, we can reformulate this as
\begin{align}
	T_m^\alpha = \sum_{l=2}^\alpha c_{\alpha l} l^m 
                     - \sum_{k,l} c_{mk} c_{\alpha l} l 
                     + \sum_{k,l} c_{mk} c_{\alpha l} k\,.
\end{align}
Using first Eq.~\eqref{eq:TrivialSumCmk} and then Eq.~\eqref{eq:WeightedSumCmk}
for the two last terms and substituting the sum on the right-hand side
by~$T_\alpha^m$, we have
\begin{align}
	T_m^\alpha - m = T_\alpha^m - \alpha\,. \label{eq:Exchange}
\end{align}
The nice property of Eq.~\eqref{eq:Exchange} is that we
always can use it to shorten or lengthen~$T_m^\alpha$ by interchanging the
indices accordingly. In particular, if $\alpha = 2$ or $\alpha = 3$ the sums
can be shortened to have a single term:
\begin{align}
	& T_m^2 - m = \sum_{k=2}^m c_{mk} k (k - 1) = 2^m - 2\,, \label{eq:MersenneConnection}\\
\intertext{or}
	& T_m^3 - m = \sum_{k=2}^m c_{mk} (k-1) k (k+1) = 3^m - 3\,.
\end{align}
We want to point out, that these equations relate the coefficients $c_{mk}$ to
the Mersenne numbers\cite{oeis2022mersenne} and to the OEIS entry
A024023\cite{oeis2022sloane}.
Also, note that $T_m^2 -m$ is coincidentally the double of $\sum_{k=2}^m c_{mk}
\binom{k}{2}$ as used by Cereceda.\cite{cereceda2022powersums} A similar
relation can be found for the case $m=3$ (cf. the identities following
Eq.~(8) in Ref.~\cite{cereceda2022powersums}).

\section{Another Identity using the closed form of the $a_{mk}$ Coefficients}
We can start from Eq.~\eqref{eq:PsiExpansion} and use
Corollary~\ref{cor:Commutator} to find
\begin{align}
	n^m = \sum_{k=2}^m c_{mk} \psi_n(k) - m + n\,.
\end{align}
On the right-hand side, we expand $\psi_n(k)$ using Eq.~\eqref{eq:nExpansion}.
Shuffling all other terms to the left-hand side, we have
\begin{align}
	n^m - n + m = \sum_{k=2}^m c_{mk} \sum_{l=2}^n a_{nl} k^l\,.
\end{align}
We switch the summation order and substitute with~$T_m^l$:
\begin{align}
	n^m - n + m = \sum_{l=2}^n a_{nl} T_m^l\,. \label{eq:PowerExpansion}
\end{align}
In particular, for the simple cases with only one non-zero $c_{mk}$,
$m=2$ and $m=3$, we have the two equations
\begin{align*}
	n^2 - n + 2 = \sum_{l=2}^n 2^l a_{nl}
	\quad{\mathrm{and}}\quad
	n^3 - n + 3 = \sum_{l=2}^n 3^l a_{nl}\,,
\end{align*}
the forms of which are strongly reminiscent of the binary and ternary
representations of a number (mind however, that the $a_{nl}$ have rational
values in general).

\section{Remark on Pythagoras' equation}
Given Eq.~\eqref{eq:PowerExpansion}, we can reformulate equations of the form
\begin{align}
	\sum_{i=1}^{N-1} q_i^m - q_N^m = 0\,,
\end{align}
insofar we have
\begin{align}
	\sum_{i=1}^{N-1} \left( \sum_{l=2}^{q_i} a_{q_i,l} T_m^l \right) 
	  - \sum_{j=2}^{q_N} a_{q_N,j} T_m^j 
	  - \sum_{i=1}^{N-1} q_i + q_N - (N-2) m = 0\,.
\end{align}

We want to consider the example $N=3$, with $q_1 = A$, $q_2 = B$, $q_3 = C$,
and $m=2$, i.e. $A^2 + B^2 - C^2 = 0$. Pythagoras' equation has then the
form\footnote{ Another representation of this equation in terms of the Mersenne
numbers comes to mind immediately: $\sum_{l=2}^{\max(A,B,C)} (2^{l-1} -1)
(a_{Al} + a_{Bl} - a_{Cl}) = \frac{C - A - B}{2}\,.$}
\begin{align}
    \sum_{l=2}^{\max(A,B,C)} 2^l (a_{Al} + a_{Bl} - a_{Cl}) &= 2 + C - A - B\,.
\end{align}
The solutions $(A,B,C)$ of this equation are---of course---the Pythagorean triples. 
It is intriguing how the simple addition and subtraction of squares is
transformed into a sum over weighted hypercubes of common sidelength and ``a
bit'', $2+C-A-B$. 

\printbibliography
\end{document}